\documentclass[microtype]{gtpart}
\agtart
\usepackage{amsmath, amssymb, amscd}

\usepackage{graphicx}
\usepackage[all,import,rotate]{xy}
\usepackage{color}
\usepackage{colortbl}

\usepackage[margin=5pt,font=small,labelfont=bf, labelsep=period]{caption}
\usepackage[margin=0pt]{subfig}

\usepackage{ragged2e}

\title{Closed surface bundles of least volume}

\author{John W.~Aaber}
\givenname{John}
\surname{Aaber}
\address{}

\email{tophat1234@hotmail}

\author{Nathan M. Dunfield}
\givenname{Nathan}
\surname{Dunfield}
\address{ Dept. of Math., MC-382 \\
          University of Illinois \\
          1409 W. Green Street \\
          Urbana, IL 61801, USA
}
\email{nathan@dunfield.info}
\urladdr{http://dunfield.info}

\keyword{surface bundle}
\keyword{pseudo-Anosov}
\keyword{minimal volume}
\keyword{minimal dilatation}

\subject{primary}{msc2010}{57M50}
\subject{secondary}{msc2010}{37E30, 37E40}

\arxivreference{arXiv:1002.3423}
\arxivpassword{h7842}

\def\RCS$#1: #2 ${\expandafter\def\csname RCS#1\endcsname{#2}}
\RCS$Revision: 80 $
\RCS$Date: 2010-11-21 15:48:06 -0600 (Sun, 21 Nov 2010) $

\renewcommand{\C}{{\mathbb C}}
\renewcommand{\Q}{{\mathbb Q}}
\newcommand{\F}{{\mathbb F}}
\renewcommand{\R}{{\mathbb R}}
\renewcommand{\Z}{{\mathbb Z}}
\renewcommand{\H}{{\mathbb H}}
\newcommand{\maps}{\colon\thinspace}
\DeclareMathOperator{\rank}{rank}
\DeclareMathOperator{\Out}{Out}
\DeclareMathOperator{\vol}{Vol}
\newcommand{\abs}[1]{{\left| #1 \right|}}
\newcommand{\norm}[1]{{\|  #1 \|}}
\newcommand{\pair}[1]{\left\langle #1 \right\rangle}

\newcommand{\spandef}[2]{{  \left\langle  {#1}  \ \left| \   {#2} \right. \right\rangle }}
\newcommand{\setdef}[2]{{  \left\{  {#1}  \ \left| \   {#2} \right. \right\} }}

\newcommand{\mtext}[1]{\quad\mbox{#1}\quad}

\newcommand{\cF}{\mathcal{F}}
\newcommand{\tcF}{\widetilde{\mathcal{F}}}

\renewcommand{\u}{\mathbf{u}}
\newcommand{\ubar}{\overline{u}}
\newcommand{\vbar}{\overline{v}}
\newcommand{\A}{\mathcal{A}}
\newcommand{\E}{\mathcal{E}}

\renewcommand{\v}{\widetilde{v}}
\newcommand{\e}{\widetilde{e}}

\newcommand{\orient}{$\circ$}
\newcommand{\htabspace}{{\hspace{0.5cm} }}

\newcommand{\mysmallmatrix}[4]{ \bigl( \begin{smallmatrix}  #1&#2\\ #3&#4 \end{smallmatrix} \bigr)
}

\theoremstyle{plain}

\swapnumbers
\newtheorem{theorem}{Theorem}[section]
\newtheorem{conjecture}[theorem]{Conjecture}
\newtheorem{lemma}[theorem]{Lemma}
\newtheorem{corollary}[theorem]{Corollary}
\newtheorem{proposition}[theorem]{Proposition}

\newtheorem{question}[theorem]{Question}

\theoremstyle{definition}
\newtheorem{definition}[theorem]{Definition}

\theoremstyle{remark}
\newtheorem{remark}[theorem]{Remark}

\newenvironment{proofsteplist}{\
  \begin{list}{}{\
      \setlength{\leftmargin}{\parindent}
      \setlength{\itemsep}{0.15cm}
    }}
  {\end{list}}

\makeatletter
  \let\c@theorem=\c@subsection
  \let\c@figure=\c@subsection
  \let\p@figure=\p@subsection
  
  \let\cl@figure=\cl@subsection
  \let\c@table=\c@subsection
  \let\p@table=\p@subsection
  
  \let\cl@table=\cl@subsection
  \let\c@equation=\c@subsection
  \let\p@equation=\p@subsection
  
  \let\cl@equation=\cl@subsection
\makeatother

\newcommand{\hyp}{\nobreakdash-\hspace{0pt}}
\newcommand{\3}[1]{3\hyp}

\newenvironment{xyoverpic}[3]
{%
\begin{xy}
\xyimport#1{\includegraphics[#2]{#3}}
}{\end{xy}}

\newenvironment{cxyoverpic}[3]
{%
\centering \leavevmode
\begin{xyoverpic}{#1}{#2}{#3}
}{\end{xyoverpic}}

\begin{document}

\begin{abstract}
  Since the set of volumes of hyperbolic 3-manifolds is well ordered,
  for each fixed $g$ there is a genus-$g$ surface bundle over the
  circle of minimal volume.  Here, we introduce an explicit family of
  genus-$g$ bundles which we conjecture are the unique such manifolds
  of minimal volume.  Conditional on a very plausible assumption, we
  prove that this is indeed the case when $g$ is large.  The proof
  combines a soft geometric limit argument with a detailed
  Neumann-Zagier asymptotic formula for the volumes of Dehn fillings.

  Our examples are all Dehn fillings on the sibling of the Whitehead
  manifold, and we also analyze the dilatations of all closed surface
  bundles obtained in this way, identifying those with minimal
  dilatation.  This gives new families of pseudo-Anosovs with low
  dilatation, including a genus 7 example which minimizes
  dilatation among all those with orientable invariant foliations.
  \end{abstract}

\maketitle

\section{Introduction}
\label{sec-intro}

\subsection{Volumes of bundles}\label{subsec-vol-bundles}
Thurston and Jorgensen proved that the set of volumes of hyperbolic
3-manifolds is a well-ordered set of type $\omega^\omega$.  In particular, any
topological class of hyperbolic 3-manifolds contains one of least
volume.  Recent work of Gabai, Meyerhoff, and Milley has
provably identified the hyperbolic manifold of least volume, and, for
instance, the smallest orientable manifolds with one and two cusps are
also known (see the recent survey \cite{GabaiMeyerhoffMileySurvey} for
more).

Here, we focus on closed hyperbolic 3-manifolds which are surface
bundles over the circle.  Fix a genus $g$, and let $\Sigma_g$ denote the
closed surface of that genus.  Consider all hyperbolic 3-manifolds
which are $\Sigma_g$-bundles over the circle.  As discussed, there is at
least one such manifold of least volume, and one goal here is to
identify these manifolds.

To this end, in Section~\ref{sec-examples} we define for each $g > 2$
a certain $\Sigma_g$-bundle $M_g$ as an explicit Dehn filling on the
Whitehead sibling manifold $W$.  For instance, if $g \equiv 0, 1 \bmod 5$,
then $M_g = W\left( \frac{g}{2g+5}, - \frac{2g+3}{g+4}\right)$.  We
conjecture that the $M_g$ are the minimal volume genus-$g$ bundles for
all $g > 2$.  For large $g$, we prove this contingent on a very
plausible conjecture about closed hyperbolic 3-manifolds with $b_1(N)
= \dim H_1(N; \Q) \geq 2$.  In particular, we show
\begin{theorem}\label{thm-main-vol}
  Assume that any closed hyperbolic 3-manifold $N$ with $b_1(N) \geq 2$
  has volume at least $V_8 \approx 3.663862$.  Then for all large $g$, the
  manifold $M_g$ is the unique $\Sigma_g$-bundle of minimal volume.
\end{theorem}
We give evidence for the needed assumption in Section~\ref{sec-qual-thm}.
Alternatively, our proof of Theorem~\ref{thm-main-vol} shows
\begin{theorem}
  Either
  \begin{enumerate}
  \item $M_g$ is the unique $\Sigma_g$-bundle of minimal volume for all
    large $g$, or
  \item There is a single closed hyperbolic \3-manifold $N$ which is
    the minimal volume $\Sigma_g$-bundle for infinitely many $g$.  (In particular, $b_1(N) \geq 2$.)
  \end{enumerate}
\end{theorem}

The proof of Theorem~\ref{thm-main-vol} has two main parts.  The
first, Theorem~\ref{thm-qual-vol}, is the qualitative statement that
the minimal volume $\Sigma_g$-bundles, for large $g$, must be Dehn
fillings on the Whitehead sibling $W$.  This is proved by combining a
soft geometric limit argument with work of Agol \cite{AgolTwoCusped}
on volumes of cusped manifolds.  Then, in
Theorem~\ref{thm-Mg-min-vol}, we sift through the large number of
$\Sigma_g$-bundles that arise from $W$ and find the one with least volume.
To do this, in Section~\ref{sec-vol} we work out the first difficult
term of the Neumann-Zagier asymptotic expansion for volumes of Dehn
fillings of $W$.

\subsection{Dilatations of pseudo-Anosovs}

We now discuss another minimization problem that turns out to be
closely related to finding the smallest volume surface bundles.  A
pseudo-Anosov homeomorphism $\psi$ of $\Sigma_g$ has an associated
dilatation $\lambda > 1$, also called the stretch-factor or growth, which
measures how the invariant foliations are distorted.  The number $\lambda$
is an algebraic integer, and $\log(\lambda)$ is the topological entropy of
the map $\psi$.  For each fixed genus, there is a smallest such
dilatation, which we denote $\delta_g$.  There has been much interest in
determining $\delta_g$, in part because $\log \delta_g$ is the length of the
shortest geodesic in the Teichm\"uller metric on the moduli space of
curves.

For a pseudo-Anosov $\psi$ of $\Sigma_g$, its mapping torus $M_\psi$ is a
hyperbolic 3-manifold.  Conversely, a hyperbolic 3-manifold which is a
$\Sigma_g$-bundle over the circle must have pseudo-Anosov monodromy.  It is
natural to ask whether there is a connection between the dilatation of
$\psi$ and the volume of $M_\psi$.   Recently, Farb, Leininger, and Margalit
have shown that, qualitatively, pseudo-Anosovs with low
dilatations come from manifolds with low volume.
\begin{theorem}[\cite{FarbLeiningerMargalit2009}]
  For each $L > 0$, there exist finitely many finite-volume hyperbolic
  3-manifolds $N_1,\ldots,N_n$ so that the following holds.  If $\psi$ is a
  pseudo-Anosov homeomorphism of a surface $\Sigma$ with $\chi(\Sigma) < 0$ and
  $\log( \lambda(\psi) ) < \frac{1}{\abs{\chi(\Sigma)}} L$, the bundle $M_\psi$
  is a Dehn filling on one of the $N_k$.  In particular, there is a $V
  = V(L)$ so that $\vol(M_\psi) \leq V$.
\end{theorem}
Moreover, they show that the fibration of $M_\psi$ comes from a
fibration of $N_k$ that extends over the solid tori added by Dehn
filling.  Thus, in particular, the bundles coming from the minimal
dilatation pseudo-Anosovs are all obtained by this \emph{fibered Dehn
  filling} from finitely many cusped hyperbolic 3-manifolds.  Here, we
focus on the smallest cusped manifold which has fibered Dehn fillings
which are $\Sigma_g$-bundles with $g$ arbitrarily large, namely the
Whitehead sibling.  This is thus a very natural place to look for
pseudo-Anosovs realizing the minimal dilatation $\delta_g$.

However, the idea of searching for pseudo-Anosovs with low dilatation
among the monodromies of low-volume manifolds predates
\cite{FarbLeiningerMargalit2009}.  In particular, earlier Kin, Kojima,
and Takasawa had asked (Question 5.1 of \cite{KinKojimaTakasawa2009}):
\begin{question}\label{quest-simult-min}
  For a fixed surface $\Sigma$, possibly with punctures, is there a
  pseudo-Anosov homeomorphism $\psi$ of $\Sigma$ which simultaneously
  minimizes $\lambda(\psi)$ and the volume of the mapping torus $M_\psi$?
\end{question}
There is good evidence that the answer is yes when $\Sigma$ is a punctured
sphere \cite{KinKojimaTakasawa2009, VenzkeThesis}.  However, by
combining our examples here with those of E.~Hironaka
\cite{Hironaka2009}, we show in
Theorem~\ref{thm-min-vol-dilatation-differ} that the answer is no for
large $g \equiv 3 \bmod 5$, provided Conjecture~\ref{conj-closed-b12}
holds.  However, we also show
\begin{theorem}
  For large $g$, the monodromy of $M_g$ minimizes dilatation among all
  $\Sigma_g$-bundles obtained from filling $W$.
\end{theorem}
\noindent
The proof uses McMullen's Teichm\"uller polynomial from
\cite{McMullen2000} to systematically calculate and analyze the
dilatation of all fillings of $W$.  Independently, Kin and Takasawa
\cite{KinTakasawa2010} have also analyzed these fillings, obtaining
similar results about dilatations (though not volumes) to those we
describe here.

Let $\lambda_g$ be the dilatation of the monodromy of $M_g$ from
Section~\ref{subsec-vol-bundles}.  For $g \equiv 3 \bmod 5$, we show that
$\lambda_g > \delta_g$.  However, for $g \not\equiv 3 \bmod 5$ and $g \neq 4$, the
dilatations $\lambda_g$ seem to be the best known and in some cases better
than previous examples.  See Table~\ref{table-dilatation} for the
known examples in low genus.  Very intriguingly, the dilatations of
both our examples and those of \cite{Hironaka2009} come from the same
family of polynomials introduced by Lanneau-Thiffeault
\cite{LanneauThiffeault2009}:
\begin{equation}\label{eq-LT}
LT_{a,b} (t) = t^{2a} - t^{a + b} - t^a - t^{a - b} + 1 \mtext{with $a > b$.}
\end{equation}
Both here and in \cite{Hironaka2009}, the parameter $a$ is within $2$
of $g$ and $b$ is small.

\begin{table}
  \definecolor{agray}{gray}{0.8}
  \begin{tabular}{rrr}

$g$ & $\lambda_g$ \htabspace & $\epsilon_g$ \htabspace \\
\cline{1-3}\\[-12.9pt]
3 & \orient  1.5061357 & \orient \cellcolor{agray} 1.4012684\\
4 &   1.2806382 &  \cellcolor{agray} 1.2612310\\
5 &  \cellcolor{agray} 1.1487947 & \orient  1.1762808\\
6 &  \cellcolor{agray} 1.1287609 &   1.1617044\\
7 & \orient \cellcolor{agray} 1.1154811 & \orient  1.1369400\\
8 &   1.1369400 &  \cellcolor{agray} 1.1135007\\
9 & \orient \cellcolor{agray} 1.0928247 & \orient  1.1054009\\
10 &  \cellcolor{agray} 1.0837668 &   1.0946555\\
11 &  \cellcolor{agray} 1.0770457 & \orient  1.0837668\\
12 &  \cellcolor{agray} 1.0726646 &   1.0787378\\
13 & \orient  1.0776673 & \orient \cellcolor{agray} 1.0726646\\
14 &  \cellcolor{agray} 1.0629875 &   1.0663959\\
15 &  \cellcolor{agray} 1.0583358 & \orient  1.0629875\\
16 &  \cellcolor{agray} 1.0549983 &   1.0590666\\
17 & \orient \cellcolor{agray} 1.0522148 & \orient  1.0549983\\
18 &   1.0568996 &  \cellcolor{agray} 1.0525403\\
19 & \orient \cellcolor{agray} 1.0470846 & \orient  1.0497943\\
20 &  \cellcolor{agray} 1.0447582 &   1.0469430\\

\end{tabular} \hspace{2cm}
\begin{tabular}{rrr}
$g$ & $\lambda_g$ \htabspace & $\epsilon_g$  \htabspace \\
\cline{1-3}\\[-12.9pt]
21 &  \cellcolor{agray} 1.0427679 & \orient  1.0450879\\
22 &  \cellcolor{agray} 1.0425372 &   1.0430555\\
23 & \orient  1.0428751 & \orient \cellcolor{agray} 1.0409471\\
24 &  \cellcolor{agray} 1.0379320 &   1.0394979\\
25 &  \cellcolor{agray} 1.0363101 & \orient  1.0379320\\
26 &  \cellcolor{agray} 1.0349894 &   1.0363101\\
27 & \orient \cellcolor{agray} 1.0338143 & \orient  1.0351470\\
28 &   1.0351470 &  \cellcolor{agray} 1.0339030\\
29 & \orient \cellcolor{agray} 1.0315904 & \orient  1.0326168\\
30 &  \cellcolor{agray} 1.0305456 &   1.0316628\\
31 &  \cellcolor{agray} 1.0296057 & \orient  1.0306505\\
32 &  \cellcolor{agray} 1.0288091 &   1.0296057\\
33 & \orient  1.0296414 & \orient \cellcolor{agray} 1.0288091\\
34 &  \cellcolor{agray} 1.0275530 &   1.0279692\\
35 &  \cellcolor{agray} 1.0263611 & \orient  1.0271037\\
36 &  \cellcolor{agray} 1.0256581 &   1.0264286\\
37 & \orient \cellcolor{agray} 1.0250132 & \orient  1.0257204\\
38 &   1.0257204 &  \cellcolor{agray} 1.0249917\\
  \end{tabular}

\vspace{0.5cm}
\caption{
  Low-dilatation pseudo-Anosovs.  Here $\lambda_g$ is the dilatation
  associated to the manifold $M_g$, and $\epsilon_g$ is the dilatation of the
  examples in Theorem~1.3 of \cite{Hironaka2009}. The darkened boxes
  have lower numbers, and $\circ$ indicates the invariant foliations are orientable.
}\label{table-dilatation}
\end{table}

\begin{table}
  \definecolor{agray}{gray}{0.8}
  \begin{tabular}{rrr}
$g$ & $\lambda_g^+$ \htabspace & $\epsilon_g^+$ \htabspace \\
\cline{1-3}\\[-12.9pt]
3 &  1.5061357 & \cellcolor{agray} $*$ 1.4012684\\
5 & $*$ 1.1762808 & \cellcolor{agray} $*$ 1.1762808\\
7 & \cellcolor{agray} $*$ 1.1154811 &  1.1369400\\
9 & \cellcolor{agray} 1.0928247 &  1.1054009\\
11 & \cellcolor{agray} 1.0803146 &  1.0837668\\
13 &  1.0776673 & \cellcolor{agray} 1.0726646\\
15 & \cellcolor{agray} 1.0597294 &  1.0629875\\
17 & \cellcolor{agray} 1.0522148 &  1.0549983\\
19 & \cellcolor{agray} 1.0470846 &  1.0497943\\
21 & \cellcolor{agray} 1.0433121 &  1.0450879\\
23 &  1.0428751 & \cellcolor{agray} 1.0409471\\
\end{tabular} \hspace{2cm}
\begin{tabular}{rrr}
$g$ & $\lambda_g^+$ \htabspace & $\epsilon_g^+$ \htabspace \\
\cline{1-3}\\[-12.9pt]
25 & \cellcolor{agray} 1.0366424 &  1.0379320\\
27 & \cellcolor{agray} 1.0338143 &  1.0351470\\
29 & \cellcolor{agray} 1.0315904 &  1.0326168\\
31 & \cellcolor{agray} 1.0297856 &  1.0306505\\
33 &  1.0296414 & \cellcolor{agray} 1.0288091\\
35 & \cellcolor{agray} 1.0264880 &  1.0271037\\
37 & \cellcolor{agray} 1.0250132 &  1.0257204\\
39 & \cellcolor{agray} 1.0237761 &  1.0244123\\
41 & \cellcolor{agray} 1.0227200 &  1.0231851\\
43 &  1.0226556 & \cellcolor{agray} 1.0221591\\
45 & \cellcolor{agray} 1.0207537 &  1.0211819\\
  \end{tabular}

\vspace{0.5cm}
\caption{
  Pseudo-Anosovs with orientable invariant foliations.
  Here $\lambda_g^+$ is the least dilatation of any example with orientable
  foliations obtained from $W$, and $\epsilon_g^+$ is that of the examples in
  Theorem~1.4 of \cite{Hironaka2009}.  The manifold
  $W$ does not generate orientable examples when $g$ is even; however
  \cite{Hironaka2009} provides such whenever $g \not \equiv 0 \bmod
  6$.  The darkened boxes have lower numbers, and a $*$ indicates
  the dilatation is known to be equal to $\delta_g^+$.
}\label{table-dilatation2}
\end{table}

However, asymptotically the $\lambda_g$ are no better than the previously
known examples.  A basic result of Penner says that $\log \delta_g \asymp
\frac{1}{g}$ \cite{Penner1991}.  In that context,  we show
\begin{theorem}
  We have  $\lim_{g \to \infty} \lambda_{g}^{g} = \frac{3+\sqrt{5}}{2} = 1 + \gamma$,
  where $\gamma$ is the golden ratio.
\end{theorem}
\noindent
The same result holds for the examples of \cite{Hironaka2009} and also
those of \cite[Thm.~1.5]{KinTakasawa2010}.

\subsection{Orientable foliations}

A natural subclass of pseudo-Anosovs is the subclass whose
invariant foliations are orientable.  Analogously, define $\delta_g^+$ to
be the least dilatation among such maps.  While $\delta_g$ is known only
when $g \leq 2$ \cite{ChoHam2008}, the work of Lanneau-Thiffeault
\cite{LanneauThiffeault2009} and Hironaka \cite{Hironaka2009} combine
to determine $\delta_g^+$ when $g \in \{2, 3, 4, 5, 8\}$.  We add one more,
showing
\begin{theorem}\label{thm-delta-plus-7}
  The monodromy $\psi_7$ of the manifold $M_7$ minimizes dilatation among
  all pseudo-Anosovs of $\Sigma_7$ with orientable invariant foliations.
  In particular, one has  $\delta_7^+ \approx 1.11548110945659$.
\end{theorem}
\noindent
Here the needed lower bound comes directly from
\cite{LanneauThiffeault2009}, and Theorem~\ref{thm-delta-plus-7} was
independently established by \cite{KinTakasawa2010}.  See
Table~\ref{table-dilatation2} for some of the known examples with
orientable invariant foliations.

\subsection{Acknowledgments.}

We thank J. O. Button, Danny Calegari, Stefan Friedl, Toby Hall, A. J.
Hildebrand, Eriko Hironaka, Eiko Kin, Sadayoshi Kojima, Erwan Lanneau,
Chris Leininger, Curt McMullen, Shigeru Mizushima, Saul Schleimer,
Mitsuhiko Takasawa, and Jean-Luc Thiffeault for helpful conversations
and correspondence.  Dunfield was partially supported by the US NSF
via grant DMS-0707136.  The authors also thank the anonymous referee
for very helpful comments and corrections, and Yuichi Kabaya for
pointing out an error in the original version of
Figure~\ref{fig-min-surface}.  The first author John W.~Aaber was
formerly known as Rupert Venzke.

\section{The sibling of the Whitehead link complement} \label{sec-sibling}

This section gives a detailed picture of the 2-cusped hyperbolic
3-manifold that is the focus of this paper, namely the sibling of the
Whitehead link complement.  This manifold, which we denote by $W$, is
the complement of the $(-2,3,8)$-pretzel link shown in
Figure~\ref{fig-peripheral-bases}, and is manifold $m125$ in the
SnapPea census \cite{CallahanHildebrandWeeks}.  Like the Whitehead
complement, it can be built by gluing the sides of a regular ideal
octahedron in $\H^3$ (see \cite{AgolTwoCusped}); thus its volume is
\[
V_8 = 4 \times (\mbox{Catalan's
  constant}) = 4 \left( \sum_{k = 0}^\infty \frac{(-1)^k}{(2k + 1)^2} \right) \approx 3.6638623767088.
\]
SnapPea \cite{SnapPea, SnapPy} finds that its fundamental group is
\[
\pi_1(W) = \spandef{x,y}{\mathit{yxyxYxYXYXyX} = 1}  \quad \mbox{where $X = x^{-1}$ and $Y = y^{-1}$.}
\]
The defining relator actually comes from a handle decomposition of
$W$, shown in Figure~\ref{fig-handle-decomp}.
\begin{figure}[t]
\centering
\subfloat[%
A handle decomposition of $W$: think of the above pattern on the
boundary of a 3-ball $B^3$, and glue the disks together to form
a genus 2 handlebody.  The arcs then form a single closed curve, along
which you attach a 2-handle corresponding to the defining relator of
$\pi_1(W)$.  ]{%
  \begin{xyoverpic}{(227,227)}{scale=0.7}{images/heegaard}
    ,(45,182)*{x}
    ,(180,182)*{y}
    ,(45,45)*{Y}
    ,(180,45)*{X}
    ,(25,144)*+!D{1}
    ,(45,138 )*+!D{2}
    ,(68 ,147 )*!DR{3}
    ,(84 ,164 )*+!R{4}
    ,(88 ,183 )*+!R{5}
    ,( 82, 201)*+!R{6}
    ,(46,89)*+!U{1}
    ,(63,83)*+!U{2}
    ,(83,64)*+!R{3}
    ,(89,45)*+!R{4}
    ,(81,24)*+!R{5}
    ,(23,83)*+!U{6}
    ,(163,143)*+!D{1}
    ,(182,138)*+!D{2}
    ,(201,144)*+!D{3}
    ,(142,203)*+!L{4}
    ,(137,183)*+!L{5}
    ,(142,162)*+!L{6}
    ,(181,89)*+!U{1}
    ,(201,83)*+!U{2}
    ,(145,24)*+!L{3}
    ,(137,46)*+!L{4}
    ,(142,64)*+!L{5}
    ,(161,84)*+!U{6}
  \end{xyoverpic}
\label{fig-handle-decomp}} \qquad %
\subfloat[The shaded disk shown is properly embedded in $B^3$ and
has algebraic intersection number 0 with the attaching curve of the
2-handle.  Thus we can add three bands inside the 2-handle to form an
orientable surface dual to $x^* + y^*$ with $\chi = -2$.]{%
  \begin{xyoverpic}{(257,251)}{scale=0.7}{images/heegaard-surface}
    ,(60,184)*{x}
    ,(196,182)*{y}
    ,(63,46)*{Y}
    ,(197,46)*{X}
    ,(39,146)*+!D{1}
    ,(61,138)*+!D{2}
    ,(82,144)*+!D{3}
    ,(100,156)*+!DR{4}
    ,(104,182)*+!R{5}
    ,(97,201)*+!R{6}
    ,(178,144)*+!D{1}
    ,(198,138)*+!D{2}
    ,(218,145)*+!D{3}
    ,(159,202)*+!L{4}
    ,(155,182)*+!L{5}
    ,(159,162)*+!L{6}
    ,(198,89)*+!U{1}
    ,(219,82)*+!U{2}
    ,(162,24)*+!L{3}
    ,(153,46)*+!L{4}
    ,(158,63)*+!L{5}
    ,(178,84)*+!U{6}
    ,(61,89)*+!U{1}
    ,(80,83)*+!U{2}
    ,(99,63)*+!R{3}
    ,(105,46)*+!R{4}
    ,(97,24)*+!R{5}
    ,(39,83)*+!U{6}
 \end{xyoverpic}
\label{fig-min-surface}}
\caption{}
\end{figure}

Notice that the relator is very symmetric in the following sense.
Realize the free group on $\{ x, y \}$ as the fundamental group of a
wedge of two circles.  Pass to the universal abelian cover and lift
the relator as shown in Figure~\ref{fig-relator}.  From the picture,
we see that both $\{ x \mapsto y , y \mapsto x \}$ and $\{ x \mapsto x
, y \mapsto Y \}$ induce automorphisms of $\pi_1(W)$; together, they
generate a subgroup of order 8 in $\Out(\pi_1(W))$, which acts on
$\pi_1(W)^\textrm{ab} = H_1(W ; \Z) = \Z^2$ as the dihedral group of
the square.  By Mostow-Prasad rigidity, these automorphisms come from
isometries of $W$.

\begin{figure}
\subfloat[The lifted relator.]{
\begin{cxyoverpic}{(187,173)}{width=90pt}{images/word}
  ,(186,53)*++!D{x}
  ,(4,169)*+!R{y}
\end{cxyoverpic}
\label{fig-relator}}
\qquad
\subfloat[The Newton polygon of $\Delta_W$ in $H_1(W; \Z)$.]{%
\begin{cxyoverpic}{(157,154)}{width=90pt}{images/newton}
  ,(163,73)*++!D{x}
  ,(78,155)*++!L{y}
\end{cxyoverpic}
\label{fig-newton-poly}}
\qquad
\subfloat[The Thurston norm unit ball $B_T \subset H^1(W;\Z)$.]{%
\begin{cxyoverpic}{(157,154)}{width=90pt}{images/thurston}
  ,(163,73)*++!D{x^*}
  ,(78,155)*++!L{y^*}
\end{cxyoverpic}
\label{fig-thurston-norm}}
\caption{}
\end{figure}

\subsection{Thurston norm} Recall that the Thurston norm on $H^1(W;
\R)$ organizes which cohomology classes correspond to fibrations over
$S^1$ (throughout this subsection, see \cite{ThurstonNorm} for
details).  Before describing the norm in this case, we fix some
notation. As a basis for $H_1(W; \Z)$, we use the images of the
generators $\{x, y\}$ of $\pi_1(W)$ under abelianization map, denoted
by the same symbols.  However, we write the operation additively,
e.g.~$x^2 y^{-3} \mapsto 2 x - 3 y$.  We'll use the dual basis $\{ x^*, y^*
\}$ for $H^1(W ; \Z)$, and show
\begin{lemma} The unit ball $B_T$ in the Thurston norm of $W$ is the
  square shown in Figure~\ref{fig-thurston-norm}, and all four sides
  are fibered faces.  In particular, every primitive element of $H^1(W
  ; \Z)$ except for $\pm(x^* \pm y^*)$ gives rise to a fibration of $W$.
\label{lem-thurston-norm}
\end{lemma}

\begin{figure}
\begin{center}
\includegraphics[scale=0.65]{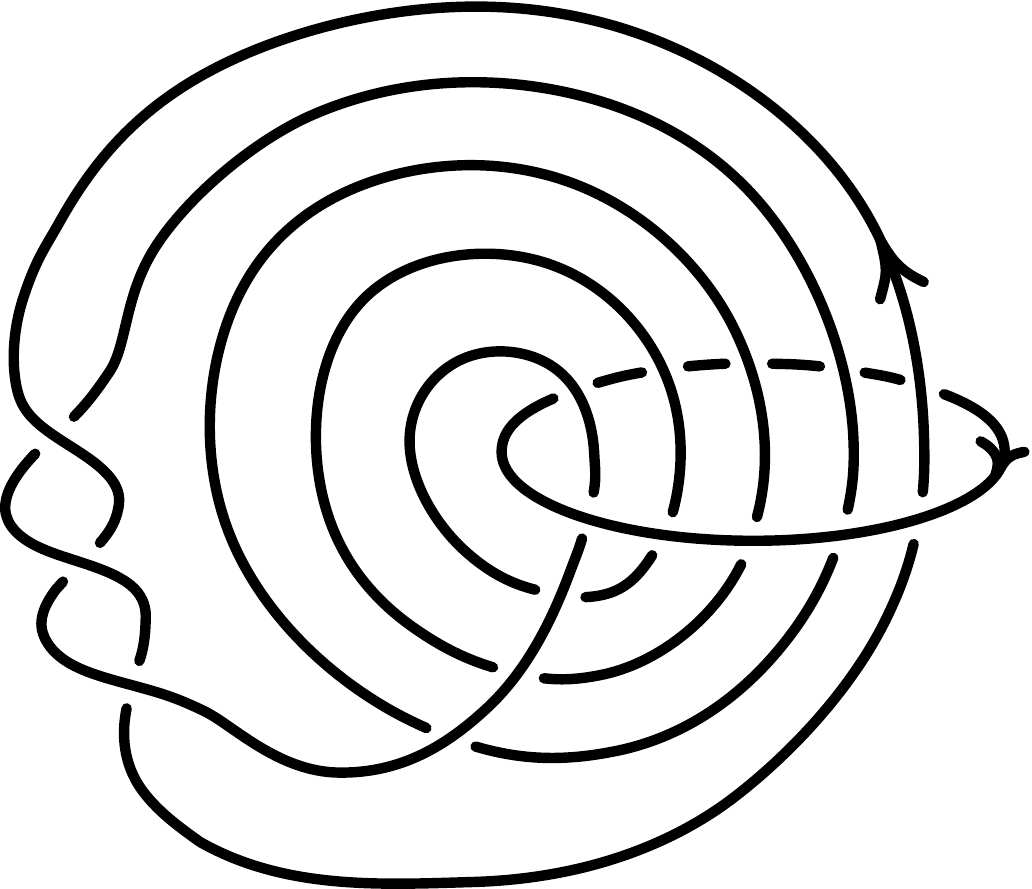}
\end{center}
\caption{The manifold $W$ is the exterior of the $(-2,3,8)$--pretzel
  link, shown here in a non-standard projection exhibiting the
  fibration of $W$ corresponding to $2 x^* + y^*$.  In terms of the
  usual homological framing on each component, we have $\mu_1 = (-7,1),
  \lambda_1 = (1,0)$ and $\mu_2 = (-4,1), \lambda_2 = (-3, 1)$, where the
  unknotted component corresponds to the second cusp.  Our orientation
  conventions are given by
  {\protect\raisebox{-2pt}{\protect\includegraphics[scale=0.35,angle=90]{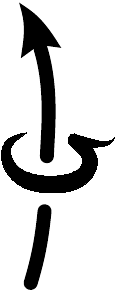}}},
  and match those of SnapPea \cite{SnapPea, SnapPy}.  That the
  fibration corresponds to $2 x^* + y^*$ follows since the meridian of the knotted component is $\lambda_1 = -x+2y$ by (\ref{eq-peripheral}).}
\label{fig-peripheral-bases}
\end{figure}

\begin{proof}
  We will use the Alexander norm introduced by McMullen
  \cite{McMullenNorm} to compute the Thurston norm.  From
  Figure~\ref{fig-relator}, it is clear the Alexander polynomial of
$W$, which lives in the group ring $\Z[ H_1(W; \Z) ]$, is
  \begin{equation}\label{eq-alex-poly}
  \Delta_W = x^{-1} + x + 1 + y^{-1} + y,
  \end{equation}
where the operation on $H_1(W; \Z)$ is written multiplicatively.  The
Newton polytope of $\Delta_W$ is shown in Figure~\ref{fig-newton-poly};
hence the unit ball $B_A$ of the Alexander norm, which is its dual, is
as in Figure~\ref{fig-thurston-norm}.  Since $b_1(W) > 1$, we have
that $B_T \subset B_A$ \cite{McMullenNorm}.  To complete the proof we need
to show that $B_T = B_A$ and all faces (sides) of $B_T$ are fibered.
By the convexity of $B_T$ and the symmetries of $W$ noted above, it is
enough to show:
\begin{enumerate}
\item The Thurston norm $\norm{x^* + y^*} \leq 2$.
\item The manifold $W$ fibers over the circle.
\end{enumerate}
For (1), we simply exhibit a twice-punctured torus dual to $x^* + y^*$
in Figure~\ref{fig-min-surface}.  For (2),
Figure~\ref{fig-peripheral-bases} shows that $W$ fibers over the
circle with fiber a 5-punctured disc, with corresponding monodromy $\psi
= \sigma_1^3\sigma_2\sigma_3\sigma_4$.
\end{proof}

\subsection{Dehn filling} \label{sec-dehn-filling}

Think of $W$ as a compact manifold having two torus boundary components
$\partial_1W$ and $\partial_2W$.  On each boundary component, we use the basis
$(\mu_k, \lambda_k)$ of $H_1(\partial_kW;\Z)$ defined in
Figure~\ref{fig-peripheral-bases}.  The shape of each cusp in the
hyperbolic structure on $\mathrm{int}(W)$ with respect to $( \mu_k,
\lambda_k)$ is a Euclidean square. For either component, we have that the
image of $H_1(\partial_kW ; \Z) \to H_1(W ; \Z) \cong \Z^2$ is a sublattice
of index $5$; explicitly, SnapPea finds that
\begin{equation}\label{eq-peripheral}
 (\mu_1, \lambda_1) \mapsto (2 x + y, -x + 2 y) \mtext{and}  (\mu_2, \lambda_2) \mapsto (x + 2 y, 2 x -  y).
\end{equation}
\begin{remark}
  This may seem like a minor point, but the fact $H_1(\partial_kW ; \Z)$ is
  nearly all of $H_1(W ; \Z)$ is what distinguishes this case from
  that of the Whitehead complement, where the image is merely $\Z$.
\end{remark}

Now given any $\phi \in H^1(W ; \Z)$, there is a unique simple closed
curve in the kernel of $\phi \maps H_1(\partial_k W) \to \Z$.  Thus, there is a
unique Dehn filling on $W$ yielding a closed manifold $W_\phi$ so that
$\phi$ extends to a map $H_1(W_\phi) \to \Z$.  Dually, if $F_\phi$ is a
surface dual to $\phi$, then $W_\phi$ is the Dehn filling of $W$ where
the components of $\partial F$ bound discs in the attached solid tori.

We are interested in those primitive $\phi \in H^1(W;\Z)$ corresponding
to fibrations.  Because of the symmetries of $W$ noted above, we will
restrict to $\phi = a x^* + b y^*$ where $0 \leq b < a$ and $\gcd(a,b) =
1$.  Because $H_1(\partial_kW;\Z)$ has index 5 in $H_1(W;\Z)$, there are
various small repeating issues throughout this paper, requiring us to
define
\[
\epsilon(\phi) =  \begin{cases} 1 & \mbox{if $a^2 + b^2 \equiv 0 \bmod 5$} \\ 0 & \mbox{otherwise.} \end{cases}
\]
\[
\epsilon_1(\phi) = \begin{cases} 1 & \mbox{if $2 a + b \equiv 0 \bmod 5$} \\ 0 & \mbox{otherwise.} \end{cases} \mtext{and} \epsilon_2(\phi) = \begin{cases} 1 & \mbox{if $a + 2 b \equiv 0 \bmod 5$} \\ 0 & \mbox{otherwise.} \end{cases}
\]
These quantities are related by

\begin{lemma}\label{lem-epsilon-props}
  Suppose $\gcd(a,b) = 1$.  Then $\epsilon(\phi) = \epsilon_1(\phi) + \epsilon_2(\phi)$.
\end{lemma}

\begin{proof}
  First notice that one never has $\epsilon_1(\phi) = \epsilon_2(\phi) = 1$, as then
 $a \equiv b \equiv 0 \bmod 5$, since $\mysmallmatrix 2 1 1
  2$ is invertible over $\F_5$.  The result then follows from the
  observation that $(2a + b)(a + 2b) \equiv 2(a^2 + b^2) \bmod 5$.
\end{proof}

\begin{theorem}\label{thm-genus}
  Let $\phi = a x^* + b y^*$ where $0 \leq b < a$ and $\gcd(a,b) = 1$.
  There is a unique fibration of $W_\phi$ over the circle, and the
  genus of the fiber is $a - 2 \epsilon(\phi)$.
\end{theorem}

\begin{proof}
  By Lemma~\ref{lem-thurston-norm}, the Thurston norm ball $B_T$ for
  $W$ is the square pictured in Figure~\ref{fig-thurston-norm}.  Our
  $\phi$ lies over the righthand face of $B_T$ where $x^* = 1/2$, and
  so $\norm{\phi} = 2 a$.  As all the faces of $B_T$ are fibered, the
  class $\phi$ gives a fibration of $W$ over the circle where the fiber
  $F$ has Euler characteristic $-\norm{\phi} = -2 a$.  Let $F_\phi$ be
  the corresponding closed surface in $W_\phi$.  Then $\chi(F_\phi) = \chi(F)
  + d_1 + d_2$, where $d_k$ is the number of boundary components of
  $F$ on $\partial_kW$.  An easy calculation using
  Lemma~\ref{lem-epsilon-props} shows that the claim
  \begin{equation}\label{eq-claim}
    d_k = 1 + 4 \epsilon_k(\phi)
  \end{equation}
  implies that the genus of $F_\phi$ is $a - 2 \epsilon(\phi)$.

  To complete the proof, we establish (\ref{eq-claim}), focusing on
  the case when $k = 1$ as the other is similar.  Using
  (\ref{eq-peripheral}), we have $d_1 = \gcd\left(\phi(\mu_1),
    \phi(\lambda_1)\right) = \gcd(2 a + b, -a + 2b) = \gcd(2 a + b, - 5 b)$,
  which is 5 if $2 a + b \equiv 0 \bmod 5$ and 1 otherwise.
  (Geometrically, these are the only two possibilities because
  $[H_1(W;\Z) : H_1(\partial_1W;\Z)] = 5$.)  That is, $d_1 = 1 + 4 \epsilon_1(\phi)$ as needed.
\end{proof}

For future reference, we record the precise Dehn filling
description of $W_\phi$.
\begin{lemma} \label{lem-filling-slopes}
  Let $\phi = a x^* + b y^*$ where $0 \leq b < a$ and $\gcd(a,b) = 1$.
  Then $W_\phi$ is the Dehn filling on $W$ given by $W\left(\frac{a -
      2b}{2 a + b}, \frac{-2a + b}{a + 2b}\right)$.  More precisely,
  the filling is along the homology classes
  \begin{align*}
  \gamma_1 &= 5^{-\epsilon_1(\phi)}\big( (a - 2b) \mu_1 + (2a + b) \lambda_1\big)  \\ 
   \gamma_2 &= 5^{-\epsilon_2(\phi)}\big( (- 2a + b) \mu_2 + (a + 2 b) \lambda_2\big)
  \end{align*}
\end{lemma}
\begin{proof}
  The surgery curve $\gamma_k$ is a primitive class killed by $\phi \maps
  H_1(\partial_kW; \Z) \to \Z$. Now one element in the kernel is $-\phi(\lambda_k)
  \mu_k + \phi(\mu_k) \lambda_k$, and as noted in the proof of
  Theorem~\ref{thm-genus}, we have $\gcd \left(\phi(\mu_k),
    \phi(\lambda_k)\right) = 1 + 4 \epsilon_k(\phi) = 5^{\epsilon_k(\phi)}$.  Thus
  \[
  \gamma_k = 5^{-\epsilon_k(\phi)} \big(-\phi(\lambda_k)  \mu_k + \phi(\mu_k) \lambda_k \big)
  \]
  is a primitive element of $H_1(\partial_k W; \Z)$, and expanding using
  (\ref{eq-peripheral}) gives the claimed formulae.
\end{proof}

\section{Examples with a given genus}
\label{sec-examples}

For a typical fixed genus $g$, there are many different fillings on $W$
which are $\Sigma_g$-bundles.  In this section, we organize these
manifolds and define our main examples $M_g$, which have the smallest
volume among all $\Sigma_g$-bundles arising from $W$ (at least when $g$
is large).  As in the previous section, by symmetry we consider only
those $\phi \in H^1(W; \Z)$ of the form $\phi = a x^* + b y^*$ with $0 \leq b
< a$ and $\gcd(a,b) = 1$.  We also assume that $g > 2$, as the
manifolds arising for $g \leq 2$ are not hyperbolic.  By
Theorem~\ref{thm-genus}, the $\phi$ for which $W_\phi$ is a
$\Sigma_g$-bundle fall into two classes
\begin{equation}\label{eq-Ag-Eg-defs}
\begin{split}
\A_g = \big\{  \phi = g x^* + b  y^* &  \ \big| \   \epsilon(\phi) = 0, \gcd(g, b) = 1 \big\}  \\
\E_g = \big\{  \phi = (g + 2) x^* + b  y^* &  \ \big|  \   \epsilon(\phi) = 1, \gcd(g+2, b) = 1 \big\}
\end{split}
\end{equation}
We show below that
\begin{proposition}\label{prop-basic-poss}
  For each $g > 2$,
  \begin {enumerate}
  \item   The set $\A_g$ is nonempty.
    \label{prop-part-Ag}
  \item    The set $\E_g$ is nonempty if and only if  $g \not\equiv 3 \bmod 5$ and $g > 4$.
    \label{prop-part-Eg}
  \end{enumerate}
\end{proposition}
Section~\ref{sec-vol} shows that manifolds coming from $\E_g$ typically have smaller
volume than those coming from $\A_g$, and volume increases as $b$
increases.  Hence we define
\begin{definition}\label{def-Mg}
  For each $g > 2$, if $\E_g$ is nonempty, we define $M_g$ to be
  $W_\phi$ where $\phi \in \E_g$ is the class where $b$ is smallest.
  Otherwise, $M_g$ comes from the $\phi \in \A_g$ where again $b$ is
  smallest.
\end{definition}
For large $g$, we will show in Theorem~\ref{thm-Mg-min-vol} that $M_g$
indeed minimizes volume among all $\Sigma_g$-bundles arising from $W$.  In
the rest of this section, we first prove
Proposition~\ref{prop-basic-poss}, and then refine it to give the
precise Dehn filling description of each $M_g$.

\begin{proof}[Proof Proposition~\ref{prop-basic-poss}]
  We begin with part (\ref{prop-part-Ag}).  If $g \equiv 0, 1, 4 \bmod 5$,
  then $\phi = g x^* + y^*$ is an element of $\A_g$.  So suppose instead
  that $g \equiv 2, 3 \bmod 5$.  Now $\A_3 = \{ 3 x^* + 2 y^*\}$, so we
  can assume that $g \geq 7$.  Then $g x^* + 5 y^* \in \A_g$, completing
  the proof of (\ref{prop-part-Ag}).

  For part (\ref{prop-part-Eg}), first suppose $g \equiv 3 \bmod 5$.  If
  $\phi = (g + 2) x^* + b y^*$ is in $\E_g$, the requirement that
  $\epsilon(\phi) = 1$ forces $b^2 \equiv 0 \bmod 5$.  But then $g + 2$ and $b$
  are both divisible by $5$, a contradiction.  So $\E_g$ is empty in
  this case.  If $g \equiv 0, 1 \bmod 5$, then $\phi = (g + 2) x^* + y^*$ is
  in $\E_g$, and if $g \equiv 2 \bmod 5$, we may simply take $b = g+1$.

  So we've reduced to when $g \equiv 4 \bmod 5$.  Now $\E_4$ is empty, so
  assume $g > 4$, in which case $a = g + 2 \geq 11$.  If $a$ is coprime
  to either 2 or 3, we simply use that for $b$, so assume $a = 2^c 3^d
  L$ where $c, d \geq 1$ and $L$ is coprime to 6.  If $L = 1$, we take
  $b = 7$.  Thus we can assume that $L > 1$, which forces $L \geq 7$ as
  $a$ is coprime to 5.

  When $L \equiv 1 \bmod 3$, we will select $b$ from among $\{ L+4, L+6,
  L+10, L+12 \}$, which are all less than $a \geq 6 L$ since $L \geq 7$.
  Each such number is coprime to 2, 3, and $L$ (for $L + 10$ this is
  because $5 \nmid L$), and hence to $a$.  As they take on four
  distinct values modulo 5, at least one is $\equiv 2,3 \bmod 5$ giving
  the needed $b$.  If instead $L \equiv 2 \bmod 3$, the same reasoning
  shows that we can choose $b$ from among $\{L+2, L+ 6, L + 8, L + 24\}$.
\end{proof}
A more detailed description is
\begin{proposition}\label{prop-phi-refined}
  Let $\phi = a x^* + b y^*$ be the class defining $M_g$.  Then
  \begin{enumerate}
  \item If $g \equiv 0,1 \bmod 5$, then $(a,b) = (g+2,1)$.
  \item If $g \equiv 3 \bmod 5$, then \label{prop-3-case}
    \begin{enumerate}
    \item $(a,b) = (g,2)$ when $g \equiv 3,13,23 \bmod 30$.
    \item $(a,b) = (g,3)$ when $g \equiv 8,28 \bmod 30$.
    \item $(a,b) = (g,5)$ when $g \equiv 18 \bmod 30$.
    \end{enumerate}
  \item If $g \equiv 2,4 \bmod 5$ with $g > 4$ then $a = g + 2$ and $b$ is
    the smallest number coprime to $g+2$ satisfying $b \equiv 2,3 \bmod 5$. \label{prop-comp-case}
  \item If $g = 4$, then $(a,b) = (4,1)$.
  \item For any $g$, the coefficient $b$ is $O(\log g)$ and is either
    1 or a prime.  \label{prop-bound-case}
  \end{enumerate}
\end{proposition}

\begin{proof}
  Everything is immediate from Proposition~\ref{prop-basic-poss} with
  the exception of (\ref{prop-3-case}) and (\ref{prop-bound-case}).
  For the first, assume $g \equiv 3 \bmod 5$ is the element of $\A_g$ for
  which $b$ is minimal.  Now $g$ is coprime to $5$, and so $b$ is at
  most $5$; since we need $\epsilon(\phi) = 0$, that is $g^2 + b^2 \not \equiv 0
  \bmod 5$, this means $b \in \{2, 3, 5\}$.  Which of these satisfy
  $\gcd(g, b) = 1$ can be determined by the class of $b$ modulo $30$,
  giving the statement above.

  The claim (\ref{prop-bound-case}) is clear in all cases except
  (\ref{prop-comp-case}), so assume that $g \equiv 2, 4 \bmod 5$ and $b$
  is the smallest number coprime to $a = g+2$ satisfying $b \equiv 2, 3
  \bmod 5$.  Note that $b$ is necessarily prime since, if not, it
  would have a prime factor which is $\equiv2,3 \bmod 5$ which would also
  be coprime to $a$.  The following argument, which we learned from
  A.~J.~Hildebrand, shows that $b$ is $O(\log g)$.

  Call a prime \emph{good} if it is $\equiv2, 3 \bmod 5$.  By minimality
  of $b$, every power $p^k < b$ of a good prime must divide $a$.  Thus
  \[
  a \geq \prod \setdef{ p^k }{ \mbox{$p$ a good prime and $k$ maximal with $p^k < b$. }}
  \]
  Taking logs and then applying the prime number theorem for
  arithmetic progressions, we get
  \[
  \log(a) \geq \sum_{\genfrac{}{}{0pt}{1}{p^k < b,}{\mbox{\scriptsize $p$ good}}} \log p  = \psi(b; 5, 2) + \psi(b; 5, 3) =  \frac{b}{2} + O\left(b e^{-c \sqrt{\log{b}}}\right) =  b \left(\frac{1}{2} + o(1) \right)
  \]
  where the terms in the center are Chebyshev's $\psi$-function, and $c > 0$ is
  an unknown constant (see e.g.~\cite[Cor.~2 on
  pg.~138]{Karatsuba1993}).  Hence $b \leq (2 + o(1)) \log (a)$ as needed.
\end{proof}

As in Section~\ref{sec-sibling}, let
$W\left(\frac{p_{1}}{q_{1}},\frac{p_{2}}{q_{2}}\right)$ denote the
manifold obtained by performing Dehn filling on $W$.  Combining
Proposition~\ref{prop-phi-refined} and Lemma~\ref{lem-filling-slopes},
we could just as well have defined $M_g$ for $g > 2$ via
\begin{enumerate}
\item If $g \equiv 0, 1$ mod $5$, then $M_{g} = W\left(\frac{g}{2 g+5},-\frac{2 g+3}{g+4}\right)$.
\item If $g \equiv 3$ mod $5$, then
  \begin{enumerate}
  \item $M_{g} = W\left(\frac{g-4}{2 g+2},\frac{2-2 g}{g+4}\right)$ when $g \equiv 3, 13, 23 \bmod 30$.
  \item $M_{g} = W\left(\frac{g-6}{2 g+3},\frac{3-2 g}{g+6}\right)$ when $g \equiv 8, 28 \bmod 30$.
  \item $M_{g} = W\left(\frac{g-10}{2 g+5},\frac{5-2 g}{g+10}\right)$ when $g \equiv 18 \bmod 30$.
\end{enumerate}
\item If $g \equiv 2, 4 \bmod 5$ and $g > 4$, then $M_{g} = W\left(\frac{-2 b+g+2}{b+2 g+4},\frac{b-2 (g+2)}{2 b+g+2}\right)$,  where $b$ is the smallest number coprime to $g+2$ satisfying $b \equiv 2,3 \bmod 5$.
\item $M_{4} = W\left(\frac{2}{9}, -\frac{7}{6}\right)$.
\end{enumerate}

\section{Proof of the basic theorem}
\label{sec-qual-thm}

Our claim that the $M_g$ are the smallest volume genus-$g$ bundles for large $g$ is
conditional on the following conjecture:

\begin{conjecture}\label{conj-closed-b12}
  Suppose $M$ is a closed hyperbolic 3-manifold with $b_1(M) \geq 2$.  Then $\vol(M) \geq V_8  \approx 3.663862$.
\end{conjecture}

This conjecture is very plausible for two reasons.  First, the smallest
example of this type known has volume $\geq 4.71$.  (The example is $M =
9^2_4(0,1)(0,1) = v1539(5,1) = s941(5,1)(5,1)$ and has the least volume
among all 0-surgeries on homologically split 2-component links with $<
15$ crossings, as well as the manifolds in (or obtained from) certain
standard censuses \cite{CallahanHildebrandWeeks,
  HallSchleimerCensus}.)  Second, if we broaden the context from
closed to finite-volume, experimentally it seems that the smallest
volume hyperbolic manifold with $b_1(M) \geq n$ always has cusps (see
\cite[\S10]{GabaiMeyerhoffMileySurvey}), and
Conjecture~\ref{conj-closed-b12} is known for cusped manifolds:

\begin{theorem}\label{thm-cusped-vol}
  Let $M$ be a finite-volume hyperbolic 3-manifold with at least one
  cusp. If $b_1(M) \geq 2$, then $\vol(M) \geq V_8$ with equality exactly
  when $M$ is the Whitehead complement or its sibling.
\end{theorem}

Theorem~\ref{thm-cusped-vol} is essentially proved in
\cite{AgolTwoCusped}, though we give the details below.  Understanding
the volumes of closed manifolds is fundamentally harder than that of
cusped manifolds, so Conjecture~\ref{conj-closed-b12} seems well
beyond present technology (see \cite{GabaiMeyerhoffMileySurvey} for a
survey).  It is worth noting that, by arguments similar to those
proving Theorem~\ref{thm-qual-vol}, there are at most finitely many
counterexamples to Conjecture~\ref{conj-closed-b12}.

\begin{proof}[Proof of Theorem \ref{thm-cusped-vol}]
  Suppose $M$ has at least two cusps; then the conclusion is exactly
  Theorem~3.6 of \cite{AgolTwoCusped}.  So suppose $M$ has one cusp.
  Since
  \[
  b_1(M) = \rank H^1(M; \Z) = \rank H_2(M, \partial M; \Z) \geq 2,
  \]
  we must have that $H_2(M; \Z) \neq 0$.  Let $\Sigma$ be a \emph{closed}
  incompressible surface representing a non-zero class in $H_2(M;
  \Z)$.  Since it is disjoint from the cusp, Theorem~3.5 of
  \cite{AgolTwoCusped} gives that either $\vol(M) > V_8$ or $M$ is
  obtained by identifying the sides of a regular ideal octahedron.
  The latter case gives finitely many possibilities, none of which
  have one cusp and $b_1(M) \geq 2$.  One can either check this by hand,
  or just use that such manifolds have been independently enumerated
  in \cite{CallahanHildebrandWeeks} and \cite{HeardPervovaPetronio}.
\end{proof}

We can now prove the basic version of the main result of this paper.

\begin{theorem}\label{thm-qual-vol}
  If Conjecture~\ref{conj-closed-b12} is true, then for every large
  $g$ the hyperbolic $\Sigma_g$-bundle of minimal volume is a Dehn filling
  on the sibling $W$ of the Whitehead complement.
\end{theorem}

\begin{proof}
  Suppose to the contrary that there is an infinite sequence $(N_n,
  g_n)$, where each $N_n$ is a minimal volume $\Sigma_{g_n}$-bundle which
  is not a Dehn filling of $W$.  Here $g_n \to \infty$ as $n \to \infty$ since
  there are finitely many manifolds of any given volume and hence
  finitely many minimal volume $\Sigma_g$-bundles for each fixed $g$.
  From Section~\ref{sec-examples}, we know that for each $g > 2$ there
  exists a $\Sigma_g$-bundle $M_g$ which is a Dehn filling on $W$.  Since
  volume strictly decreases under Dehn filling, we have $\vol(M_g) <
  \vol(W) = V_8$ for all $g$, and hence by minimality that $\vol(N_n)
  < V_8$ for all $n$.  As their volumes are bounded, we pass to a
  subsequence where the $N_n$ have a geometric limit $N$; for ease of
  notation, we reindex so that $N_n$ now denotes that subsequence.  As
  $N$ is the geometric limit of the $N_n$, we have $\vol(N) = \lim_{n
    \to \infty}\vol(N_n)$ and so $\vol(N) \leq V_8$.  There are three cases,
  depending on how many cusps $N$ has.
  \begin{proofsteplist}
  \item \textbf{No cusps:} If $N$ is compact, then the sequence is
    eventually constant, with $N_n \cong N$ for all large $n$.  As $g_n
    \to \infty$, the manifold $N$ must fiber in infinitely many ways, in
    particular $b_1(N) \geq 2$.  But then
    Conjecture~\ref{conj-closed-b12} implies that $\vol(N) \geq V_8$, a
    contradiction since $\vol(N_n) < V_8$ for all $n$.

  \item \textbf{One cusp:} Suppose $N$ has one cusp.  For all but one
    Dehn filling of $N$, the first betti number of the filled manifold
    is $b_1(N) - 1$.  (The exception is the filling along the kernel
    of $H_1(\partial N; \Q) \to H_1(N; \Q)$.) Thus since $b_1(N_n) \geq 1$ for
    every $n$, we have $\nobreakspace b_1(N) \geq 2$.  But then
    Theorem~\ref{thm-cusped-vol} implies that $\vol(N) > V_8$, a
    contradiction to our observation that $\vol(N) \leq V_8$.

  \item \textbf{At least 2 cusps:} As noted, we have $\vol(N) \leq V_8$.
    Since $N$ has at least 2 cusps, we know $b_1(N) \geq 2$ and so
    Theorem~\ref{thm-cusped-vol} implies that $N$ is either the
    complement of the Whitehead link or its sibling.  Our assumption
    on the $N_n$ means it must be the former, which we denote $W'$.
    For each of the two boundary components $\partial_iW'$, there is a
    natural meridian-longitude pair $(\mu_i, \lambda_i)$ in $H_1(\partial_iW' ;
    \Z)$.  Here, $H_1(W' ; \Z) \cong \Z^2 = \left\langle \mu_1, \mu_2\right\rangle$
    and both $\lambda_i = 0$ in $H_1(W' ; \Z)$.  Now each $H^1(N_n; \Z) \neq
    0$, and so gives rise to a non-zero class $\phi_n \in H^1(W'; \Z)$.
    The Dehn filling curve on $\partial_1W'$ lies in the kernel of $\phi_n$,
    and thus either is $\lambda_1$ or $\phi_n = \mu_2^*$.  Since the $N_n$
    geometrically limit on $W'$, the Dehn filling curves must be
    varying on both cusps, forcing $\phi_n = \mu_2^*$ for all but
    finitely many $n$.  But the symmetric argument with the other cusp
    shows that $\phi_n = \mu_1^*$ for all but finitely many $n$, a
    contradiction.
  \end{proofsteplist}
  In each case we have a contradiction, and thus all but finitely many
  of the $N_n$ must be fillings of the Whitehead sibling $W$.
\end{proof}

\section{Volumes of Dehn fillings}
\label{sec-vol}

Since they are Dehn fillings, each $W_\phi$ has smaller volume
than $W$, and most of them have volumes very close to $\vol(W)$.  In
this section, we give very precise estimates on the volumes of most
$W_\phi$, allowing us to show
\begin{theorem}\label{thm-Mg-min-vol}
  For all large $g$, the manifold $M_g$ minimizes volume among all
  $\Sigma_g$-bundles which are Dehn fillings on the sibling of the
  Whitehead complement.
\end{theorem}
\noindent
Together with Theorem~\ref{thm-qual-vol}, this will prove
Theorem~\ref{thm-main-vol}.

Throughout this section, we consider $\phi = a x^* + b y^* \in H^1(W;\Z)$
with $0 \leq b < a$ and $\gcd(a,b) = 1$.   The key quantity controlling
$\vol(W_\phi)$ turns out to be
\[
n(\phi) = \frac{5\left(a^2 + b^2\right)}{1 + 25^{\epsilon(\phi)}}
\]
As a first step, we use work of Neumann and Zagier
\cite{NeumannZagier} to prove
\begin{lemma}\label{lem-basic-vol}
  \[
  \vol(W_\phi) = \vol(W) - \frac{\pi^2}{n(\phi)} + O\left( \frac{1}{n(\phi)^2} \right)
  \]
\end{lemma}
This lemma suggests that $n(\phi) > n(\psi)$ implies $\vol(W_\phi) >
\vol(W_\psi)$, which, if true, easily implies
Theorem~\ref{thm-Mg-min-vol} given the formula for $n(\phi)$.  However,
we will need to refine the lemma substantially in order to actually
prove the theorem.
\begin{proof}[Proof of Lemma~\ref{lem-basic-vol}]
  For any finite-volume hyperbolic 3-manifold $M$, Neumann and Zagier
  proved the following asymptotic formula for the volumes of large
  Dehn fillings.  Let $T_k$ be a Euclidean torus associated to the
  $k^\mathrm{th}$ cusp, and let $\gamma_k = p_k \mu_k + q_k \lambda_k$ be the
  Dehn filling curve for that cusp.  In terms of the quadratic form
  \[
  Q(\gamma_k) = Q(p_k, q_k) = \left(\mbox{length of $\gamma_k$}\right)^2 \big/ \big(\mbox{area of $T_k$}\big)
  \]
  Theorem 1A of \cite{NeumannZagier} states
  \begin{equation}\label{eq-NZ}
  \vol\big(\mbox{$M$ filled along the $\gamma_i$}\big) = \vol(M) - \pi^2 \sum_k \frac{1}{Q(p_k, q_k)} + O\left(\sum_k \frac{1}{p_k^4 + q_k^4}\right)
  \end{equation}
  As mentioned in Section~\ref{sec-dehn-filling}, for $W$ with
  the bases we've chosen, the $T_k$ are both squares with sides
  parallel to $\mu_k$ and $\lambda_k$.  From Lemma~\ref{lem-filling-slopes},
  one finds $Q(\gamma_k) = 5^{1 - 2 \epsilon_k(\phi)}(a^2 + b^2)$.  Using
  Lemma~\ref{lem-epsilon-props}, we see that the first sum in
  (\ref{eq-NZ}) is simply $1/n(\phi)$ as claimed.  As for the
  error term, just observe that because of the square cusp shape one
  has $p_k^4 + q_k^4 \leq (p_k^2 + q_k^2)^2 = Q(\gamma_k)^2$ and so
  \[
  \sum_k \frac{1}{p_k^4 + q_k^4} \leq \sum_k \frac{1}{Q(\gamma_k)} = \frac{1 + 5^{4 \epsilon}}{25 (a^2 + b^2)^2} \leq \frac{1}{n(\phi)^2}
  \]
  and hence we can view the error in (\ref{eq-NZ}) as $O\big(n(\phi)^{-2}\big)$.
\end{proof}

Lemma~\ref{lem-basic-vol} is enough to distinguish between the
two classes of examples introduced in (\ref{eq-Ag-Eg-defs}).

\begin{theorem}\label{thm-Ag-vs-Eg}
  For all large $g$, if $\phi \in \A_g$ and  $\psi \in \E_g$ then $\vol(W_\phi) > \vol(W_\psi)$.
\end{theorem}

\begin{proof}
  Notice that
  \[
  \frac{ n(\phi) }{ n(\psi)} = \frac{ \frac{5}{2}\left( g^2 + b_1^2
    \right)}{\frac{5}{26} \left( (g + 2)^2 + b_2^2 \right)} = 13 \frac{g^2 + b_1^2}{ (g+2)^2 + b_2^2} > \frac{13}{2} \left(\frac{g}{g + 2}\right)^2
  \]
  and thus $n(\phi) \geq 6 n(\psi)$ for large $g$.  Hence
  \[
  \vol(W_\phi) - \vol(W_\psi) = \pi^2 \left( \frac{1}{n(\psi)} - \frac{1}{n(\phi)} \right) + O\left(\frac{1}{n(\phi)^2}\right) \geq \frac{5 \pi^2}{n(\phi)} + O\left(\frac{1}{n(\phi)^2}\right)
  \]
  which is positive for large $n(\phi)$, and hence for large $g$, as desired.
\end{proof}
Unfortunately, Lemma~\ref{lem-basic-vol} is not enough to identify the
minimal volume manifold among those coming from $\A_g$ (or $\E_g$).
As a concrete example, if $g \equiv 13 \bmod 30$, then $\phi = g x^* + 2
y^*$ and $\psi = g x^* + 3 y^*$ are both in $\A_g$.  However, the
difference in volumes is governed by $1/n(\psi) - 1/n(\phi) =
2/\big(n(\phi)n(\psi)\big)$, which is comparable to the error.  To
distinguish such examples, we use that $\phi$ and $\psi$ are very closely
related (in this case, $b$ differs by $1$) to show that the error
terms in Lemma~\ref{lem-basic-vol} are also correlated, and this
difference is accurate to one higher order than the original estimate
itself.  To execute this we need the following refined volume
estimate.
\begin{theorem}\label{thm-good-vol-formula}
  Consider $\phi = a x^* + b y^* \in H^1(W;\Z)$ with $0 \leq b < a$ and $\gcd(a,b) = 1$.   If $\epsilon(\phi) = 0$ then
  \[
    \vol(W_\phi) = \vol(W) - \frac{2 \pi^2}{5(a^2 + b^2)}   - \frac{ \pi^4 \left(7a^4 + 144a^2b^2 + 7b^4\right)}{1875\left(a^2 + b^2\right)^4} +  O\left( \frac{1}{a^6} \right)
  \]
  If instead $\epsilon(\phi) = 1$ then $\vol(W_\phi) - \vol(W)$ is given by
  \[
   - \frac{26 \pi^2}{5(a^2 + b^2)}  + \frac{\pi^4\left(901(a^4 + b^4) + (-1)^{\epsilon_1(\phi)}2808(ab^3 - a^3b) - 6336a^2b^2\right)}{375\left(a^2 + b^2\right)^4} +  O\left( \frac{1}{a^6} \right)
  \]
\end{theorem}
\begin{proof}
  We will derive a more refined asymptotic expansion for $\vol(W_\phi)$
  using the method of \cite{NeumannZagier}.  For the standard
  triangulation of $W = m125$ with four tetrahedra, the gluing
  equations are
  \begin{equation}\label{eq-gluing}
     z_1(1 -z_2)(1 - z_3)(1-z_4) = -  (1 - z_1) \quad \mbox{and}
     \quad   (1 - z_1)z_2z_4 =  z_3
  \end{equation}
  where $z_i \in \C$ is the shape of the $i^{\mathrm{th}}$ tetrahedron.
  A solution to these equations gives a (typically incomplete)
  hyperbolic structure on $W$, with the induced holonomy $h$ on the
  boundary given by
  \begin{equation}\label{eq-holo}
    \begin{split}
      U_1 = h(\mu_1) =  \frac{1}{z_1 (1 - z_1)  z_2 z_4 (1 - z_4)} \quad &\mbox{and} \quad V_1 = h(-\lambda_1) = - \frac{ z_2(1-z_3)(1 - z_4)}{z_4} \\
      \\
      U_2 = h(\mu_2) =  - \frac{(1 - z_1)  z_3}{ (1-z_3)(1 - z_4)} \quad &\mbox{and} \quad V_2 = h(-\lambda_2) =  \frac{ z_1 z_2 (1-z_4)}{(1 - z_1)z_4}
      \end{split}
  \end{equation}
  Here, to match the conventions of \cite{NeumannZagier}, we are using
  bases $(\mu_k, -\lambda_k)$ for the homology of the $k^{\mathrm{th}}$ cusp.

  Following \cite{NeumannZagier}, we consider $u_k = \log U_k$ and
  $v_k = \log V_k$.  On a neighborhood of $\mathbf{0} \in \C^2$, the
  pair $\u = (u_1,u_2)$ parameterize solutions to the gluing equations
  near the complete hyperbolic structure.  The $v_k$ can be related to
  the $u_k$ via a potential function $\Phi(\u)$, where in
  particular
  \[
  v_k = \frac{1}{2} \frac{\partial\Phi}{\partial u_k}
  \]
  Because there is an isometry of $W$ interchanging $(\mu_1, \lambda_1)$ with
  $(\mu_2, \lambda_2)$, the potential function has the symmetric form
  \begin{equation}\label{eq-Phi}
    \Phi(\u) = c_1 \big(u_1^2 + u_2^2\big) + c_2 \big(u_1^4 + u_2^4\big) + c_3 u_1^2 u_2^2 + O\big( |\u|^6 \big)
  \end{equation}
The function $\Phi$ also determines the following refined volume
estimate of \cite{NeumannZagier}.  If $W_\u$ denotes the hyperbolic
structure on $W$ corresponding to $\u$, and $\Phi_k$ the terms of $\Phi$
of total degree $k$, then
\begin{equation}\label{eq-basic-vol}
\vol(W_\u) = \vol(W) + \frac{1}{4} \cdot \mathrm{Im} \left(\sum_{k=1}^2  u_k \vbar_k - \frac{1}{2} \sum_{k = 2}^\infty (k  - 2)\Phi_k(\u) \right).
\end{equation}
Thus the key is to find the coefficients $c_k$ in (\ref{eq-Phi}); the
first one is determined by the cusp shapes, which are both square, and
hence $c_1 = i$.

Using Gr\"obner bases \cite{Macaulay2} to eliminate variables from
(\ref{eq-gluing}) and (\ref{eq-holo}), we find the following
relationship between $V_1, U_1,$ and $U_2$
\begin{multline}\label{eq-horror}
U_2V_1^2  + U_1^6U_2V_1^4 - 2U_1U_2V_1^2  \left(3 - V_1 + V_1^2\right) - 2U_1^5U_2V_1^2\left(1 - V_1 + 3V_1^2\right) \\
- U_1^4\Big(\left(-1 + V_1\right)^4V_1 + U_2^2\left(-1 + V_1\right)^4V_1 +  U_2\left(-1 + 4V_1 - 15V_1^2 + 12V_1^3 - 15V_1^4\right)\Big) \\
+ U_1^2V_1\Big(-\left(-1 + V_1\right)^4 - U_2^2\left(-1 + V_1\right)^4 +  U_2V_1\left(15 - 12V_1 + 15V_1^2 - 4V_1^3 + V_1^4\right)\Big) \\
+ 2U_1^3V_1\Big(\left(-1 + V_1\right)^4 + U_2^2\left(-1 + V_1\right)^4 -  U_2\left(1 + 4V_1 + 4V_1^3 + V_1^4\right)\Big) = 0
\end{multline}
If we substitute $U_k = e^{u_k} = 1 + u_k + u_k^2/2 + u_k^3/6 + O\big(\abs{\u}^4\big)$ and
\begin{align*}
V_1 = e^{v_1}  &= \exp\left(i u_1 + 2 c_2 u_1^3 + c_3u_1u_2^2 + O\big(\abs{\u}^4\big)\right) \\
 & = 1 + i u_1 - \frac{1}{2}u_1^2 + \frac{12 c_2 -  i}{6} u_1^3 +  c_3 u_1u_2^2  + O\big(\abs{\u}^4\big)
\end{align*}
into (\ref{eq-horror}), we get
\[
-\big(1 + (8 - 8 i) c_3\big) u_1^6 u_2^2 - \left(\frac{1 + 2 i}{3} + (16 - 16 i) c_2\right)  u_1^8+ O\big(\abs{\u}^9\big)  = 0
\]
which forces $c_2 = (-3 + i)/96$, and $c_3 = -(1 + i)/16$.   In particular,
\[
v_1 = i u_1 + \frac{-3 + i}{48}u_1^3 - \frac{1 + i}{16}u_1u_2^2 +  O\big(\abs{\u}^4\big).
\]
Now (\ref{eq-basic-vol}) becomes
\begin{equation}\label{eq-big-vol}
\begin{split}
\vol(W_\u) = \vol(&W)  -  \frac{1}{4}\left( u_1 \ubar_1 + u_2 \ubar_2 \right)  -  \frac{1}{64} \mathrm{Im}\Bigg[ \frac{3+ i}{3}\big(u_1\ubar_1^3 + u_2\ubar_2^3\big)\\
  &+ (1 -  i)\big(u_1\ubar_1\ubar_2^2 + \ubar_1^2 u_2\ubar_2\big) -  \left( \frac{3 -  i}{6}\big(u_1^4 + u_2^4\big) + (1 + i)u_1^2u_2^2\right) \Bigg] \\
&\hspace{3cm}  +  O\big(\abs{\u}^6\big).
\end{split}
\end{equation}

We can now use this formula to estimate $\vol(W_\phi)$ as follows;
here, we do the case $\epsilon(\phi) = 0$, and leave the other to the reader.
By Lemma~\ref{lem-filling-slopes}, $W_\phi$ is obtained by filling along
\[
\gamma_1 =  (a - 2b) \mu_1 + (2a + b) \lambda_1 \mtext{and}   \gamma_2 =  (- 2a + b) \mu_2 + (a + 2 b) \lambda_2
\]
Thus if $\u$ is the parameter realizing the hyperbolic structure on
$W_\phi$ we have, remembering that $v_k = \log \big(h(-\lambda_k) \big)$, the
following relationship
\begin{equation}\label{eq-holo-cond}
 (a - 2b) u_1 - (2a + b) v_1 = 2 \pi i \mtext{and}  (- 2a + b) u_2 - (a + 2 b) v_2 = 2 \pi i
\end{equation}
Following \cite[\S6]{NeumannZagier}, we introduce quantities
\[
A_1 = (a - 2b)  -  (2a + b) i  \mtext{and} A_2 = (- 2a + b)  - (a + 2 b) i
\]
and use (\ref{eq-holo-cond}) to find that
\begin{equation}
\begin{split}
u_1 & = \frac{2\pi i}{A_1} -   \frac{ \pi^3(2 a + b) i}{2 A_1^2} \left(\frac{1 - i/3}{A_1^2} + \frac{1 + i}{A_2^2} \right) + O\left(\frac{1}{a^5}\right)\\
u_2 & = \frac{2\pi i}{A_2} -   \frac{ \pi^3(a + 2 b) i}{2 A_2^2} \left(\frac{1 +i}{A_1^2} + \frac{1 - i/3}{A_2^2} \right) + O\left(\frac{1}{a^5}\right)\\
\end{split}
\end{equation}
Substituting this into (\ref{eq-big-vol}) and expanding gives
  \[
    \vol(W_\phi) = \vol(W) - \frac{2 \pi^2}{5(a^2 + b^2)}   - \frac{ \pi^4 \left(7a^4 + 144a^2b^2 + 7b^4\right)}{1875\left(a^2 + b^2\right)^4} +  O\left( \frac{1}{a^6} \right)
  \]
  as claimed.
\end{proof}
We now use this last result to prove the following, which immediately
gives Theorem~\ref{thm-Mg-min-vol} when combined with
Theorem~\ref{thm-Ag-vs-Eg} and Definition~\ref{def-Mg}.
\begin{theorem}
  For all large $g$ the following is true.  Suppose $\phi = a x^* + b_1
  y^*$ and $\psi = a x^* + b_2 y^*$ are both in $\A_g$
  or both in $\E_g$.  If  $b_1 > b_2$ then $\vol(W_\phi) > \vol(W_\psi)$.
\end{theorem}
\begin{proof}
  First suppose $\phi$ and $\psi$ are in $\A_g$.  By the part of
  Theorem~\ref{thm-good-vol-formula} that followed directly from
  Lemma~\ref{lem-basic-vol}, we have
  \[
  \vol(W_\phi) - \vol(W_\psi) = \frac{2 \pi^2}{5} \frac{b_1^2 -
    b_2^2}{\left(a^2 + b_1^2\right)\left(a^2 + b_2^2 \right)} +
  O\left(\frac{1}{a^4}\right) \geq \frac{\pi^2}{10} \frac{b_1^2 -
    b_2^2}{a^4} +   O\left(\frac{1}{a^4}\right)
  \]
  Suppose the error term above is $< C/a^4$.  As $b_1 > b_2$ are
  integers, $b_1^2 - b_2^2 > 2 b_2$, so if $b_2 > 10 C $ then
  $\vol(W_\phi) > \vol(W_\psi)$.  Thus we can assume that $b_1$ and $b_2$
  are uniformly bounded, independent of $a$.

  Now the more refined part of Theorem~\ref{thm-good-vol-formula}
  gives that
  \[
    \vol(W_\phi) - \vol(W_\psi) = \frac{2 \pi^2}{5} \frac{b_1^2 - b_2^2}{\left(a^2 + b_1^2\right)\left(a^2 + b_2^2 \right)} + \frac{\pi^4 P(a, b_1, b_2)} {\left(a^2 + b_1^2\right)^4\left(a^2 + b_2^2 \right)^4} +
  O\left(\frac{1}{a^6}\right)
  \]
  where $P \in \Z[a,b_1,b_2]$ has degree $10$ in $a$.  Since the
  $b_i$ are uniformly bounded, this gives
    \[
    \vol(W_\phi) - \vol(W_\psi) = \frac{2 \pi^2}{5} \frac{b_1^2 -
      b_2^2}{\left(a^2 + b_1^2\right)\left(a^2 + b_2^2 \right)} +
    O\left(\frac{1}{a^6}\right) \geq \frac{\pi^2}{10 a^4}  + O\left(\frac{1}{a^6}\right)
  \]
  and hence $\vol(W_\phi) > \vol(W_\psi)$ for large $g$.

  If instead $\phi$ and $\psi$ are in $\E_g$, the only significant
  difference is that $P$ has degree $11$ in $a$ and thus one gets,
  when the $b_i$ are bounded, that
  \[
  \vol(W_\phi) - \vol(W_\psi) = \frac{26 \pi^2}{5} \frac{b_1^2 -
    b_2^2}{\left(a^2 + b_1^2\right)\left(a^2 + b_2^2 \right)} +
  O\left(\frac{1}{a^5}\right)
  \]
  which still suffices to prove the theorem.
\end{proof}

\section{The monodromy of  $W_\phi$}
\begin{figure}
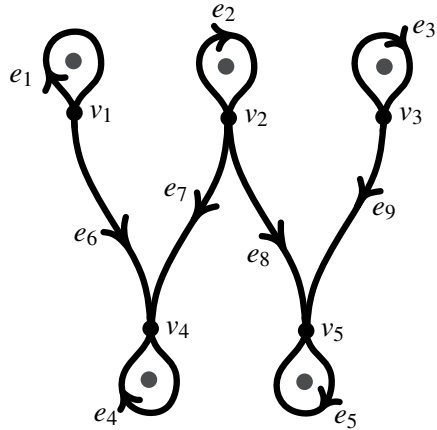

  \begin{center}
  \begin{cxyoverpic}{(175,186)}{scale=0.8}{images/traintrack}
    ,(15,143)*++!L{v_1}
    ,(88,141)*++!L{v_2}
    ,(161,142)*++!L{v_3}
    ,(51,41)*++!L{v_4}
    ,(124,40)*++!L{v_5}
    ,(4,160)*++!R{e_1}
    ,(86,179)*++!D{e_2}
    ,(171,175)*+!DL{e_3}
    ,(40,8)*+!RU{e_4}
    ,(134,7)*+!LU{e_5}
    ,(38,85)*+++!R{e_6}
    ,(73,98)*+!RD{e_7}
    ,(112,83)*+!RU{e_8}
    ,(152,105)*+!UL{e_9}
  \end{cxyoverpic}
  \end{center}
  \caption{The invariant train track for $\psi = \sigma_1^3\sigma_2\sigma_3\sigma_4$.}
  \label{fig-train-track}
\end{figure}
In this section, we compute the Teichm\"uller polynomial for a face of
the Thurston norm ball of $W$, and use this to calculate the
dilatation for each $M_g$, as well as the other closed 3-manifolds
fibering over the circle which are obtained by filling $W$.  For
background and definitions, see McMullen's original paper
\cite{McMullen2000}.   Specifically, in the notation of
Section~\ref{sec-sibling}, we have:
\begin{theorem}\label{thm-dilatation}
  Consider $\phi = a x^* + b y^* \in H^1(W;\Z)$ with $0 \leq b < a$ and
  $\gcd(a,b) = 1$.  If $(a,b) \notin \{ (1,0), (2,1), (3,1), (4,3) \}$,
  then $W_\phi$ is a hyperbolic 3-manifold which is a $\Sigma_g$-bundle,
  where $g = a - 2\epsilon(\phi)$.  The invariant foliations for the monodromy
  have $5^{\epsilon_k(\phi)}$ singularities with $5^{-\epsilon_k(\phi)} (2 a
  +(-1)^{k+1} b)$ prongs for $k = 1, 2$.  The dilatation of the
  monodromy is the largest real root of
  \[
   \Lambda_{a,b}(t) = t^{2a} - t^{a+b} - t^{a} - t^{a-b} + 1.
   \]
   Finally, the invariant foliations are (co)orientable exactly when
   $a$ is odd and $b$ is even.
\end{theorem}
Combining this with Proposition~\ref{prop-phi-refined} immediately gives:
\begin{corollary}\label{cor-dilatations}
The dilatation $\lambda_{g}$ of the monodromy of $M_{g}$ is the largest real root of $a_{g}(t)$ where
\begin{enumerate}
\item If $g \equiv 0,1 \bmod 5,$ then $a_{g}(t) = \Lambda_{g+2,1}(t) = t^{2g+4}-t^{g+3}-t^{g+2}-t^{g+1}+1$
\item If $g \equiv 3 \bmod 5,$ then
\begin{enumerate}
\item $a_{g}(t) = \Lambda_{g,2} = t^{2g} - t^{g+2} - t^{g} - t^{g-2} + 1$ when $g \equiv 3,13,23 \bmod 30$
\item $a_{g}(t) = \Lambda_{g,3} = t^{2g} - t^{g+3} - t^{g} - t^{g-3} + 1$ when $g \equiv 8,28 \bmod 30$
\item $a_{g}(t) = \Lambda_{g,5} = t^{2g} - t^{g+5} - t^{g} - t^{g-5} + 1$ when $g \equiv 18 \bmod 30$
\end{enumerate}
\item If $g \equiv 2,4 \bmod 5,$ then $a_{g}(t) = \Lambda_{g+2,b} =
  t^{2g+4} - t^{g+b+2} - t^{g+2} - t^{g-b+2} + 1,$ where $b$ is chosen
  to be the smallest number coprime to $g+2$ satisfying $b \equiv 2,3
  \bmod 5.$
\end{enumerate}
\end{corollary}

\subsection{The Teichm\"uller Polynomial}
Our starting point for Theorem~\ref{thm-dilatation} is to view $W$ as
a mapping torus of the 5-punctured disk $D$ with monodromy $\psi =
\sigma_1^3\sigma_2\sigma_3\sigma_4$ (see Figure~\ref{fig-peripheral-bases}).  Using
the Bestvina-Handel algorithm, as implemented by \cite{HallTrains}, we
found the invariant train track $\tau$ shown in
Figure~\ref{fig-train-track}.  The monodromy $\psi$ permutes the
vertices of $\tau$ via the 5-cycle $(1,3,5,4,2)$, and just permutes the
vertex loops in the compatible way.  On the other edges, $\psi$ wraps
them as follows, where ``$\cdot$'' denotes concatenation of paths from
left to righthand (as with the usual conventions for fundamental group):
\begin{equation}\label{eq-edge-action}
  e_6 \mapsto e_9 \cdot e_5 \cdot e_8^{-1} \qquad   e_7 \mapsto e_6 \cdot e_4^{-1} \cdot e_7^{-1}  \qquad e_8 \mapsto e_6 \qquad e_9 \mapsto e_8^{-1} \cdot e_2^{-1} \cdot e_7
\end{equation}

Following \cite{McMullen2000}, we can compute the Teichm\"uller
polynomial by considering the lift of $\psi$ to a certain abelian cover
$\widetilde{D} \to D$, and then considering its action on the space of
weights of $\widetilde{\tau}$.  This is analogous to computing the
Burau representation or the Alexander polynomial, and what follows can
be thought of in terms of Fox calculus.  Here, $\widetilde{D}$ comes
from the map from $\pi_1(D)$ to the infinite cyclic group $\pair{t}$
which maps a clockwise loop about a single puncture to $t$.  Observing
how the $t$-grading shifts along the image edges in
(\ref{eq-edge-action}), we find that in $\widetilde{D}$:
\[
\v_1 \mapsto \v_3 \qquad \v_2 \mapsto t^2 \v_1 \qquad \v_3 \mapsto t^3 \v_5 \qquad \v_4 \mapsto t \v_2 \qquad \v_5 \mapsto t^2 \v_4
\]
The action of $\widetilde{\psi}$ on the first five edges is determined
by the vertex action, and the others are
\begin{align*}
\e_6 &\mapsto t \e_5 + t \e_8 + \e_9 &  \e_7 &\mapsto t^2 \left(\e_4 + \e_6 + t^{-1} \e_7 \right) \\
\e_8 &\mapsto t^2 \e_6  &  \e_9 &\mapsto t^3 \left( \e_2 + t^{-1} \e_7 + \e_8 \right)
\end{align*}
Now by \cite{McMullen2000}, to get the Teichm\"uller polynomial we
simply compute
\[
\Theta(u, t) = \frac{\det \left(uI - P_{E}(t)\right)}{\det \left(uI - P_{V}(t)\right)} = -ut^5 + t^6 - u^2t^3 + u^4 - u^3t, 
\]
where $u$ is a certain element of $H_1(W)$ which goes once around the
bundle.  To put this in our usual basis for $H_1(W)$, first note from
Figure~\ref{fig-peripheral-bases} that $t = x - 2y$ and the element of
$H^1(W;\Z)$ corresponding to this fibration is $\phi = 2 x^* + y^*$.
Since $\phi(u) = 1$ and $\Theta$ determines the righthand face of the
Thurston norm ball given in Figure~\ref{fig-thurston-norm}, it follows
that $u = 2 x - 3 y$ and
\begin{equation}\label{eq-teichmuller}
\Theta(x,y) = x + x^{-1} - 1 - y - y^{-1} = -\Delta_W(-x, y).
\end{equation}

\subsection{Invariant foliations}

Consider $\phi = a x^* + b y^* \in H^1(W ; \Z)$ over the interior of the
righthand fibered face $F$, i.e. $a > 0$ and $|b| < a$.  Consider the
monodromy $\psi \maps \Sigma \to \Sigma$ of the corresponding fibration of $W$
(here $\Sigma$ is a certain surface with boundary).  We wish to understand
the topology of the invariant foliations of $\psi$ so that we can
understand the bundle structure of the Dehn-filled manifold $W_\phi$,
and use (\ref{eq-teichmuller}) to compute the stretch factor.

For concreteness, let us focus on the stable foliation $\cF \subset
\Sigma$.  This foliation can be suspended via $\psi$ to a 2-dimensional
foliation $\tcF$ of $W$ with a finite number of singular leaves.
Fried showed that $\tcF$ only depends, up to isotopy, on the
fibered face and not the particular choice of $\phi$ \cite{Fried1979}.  Thus from the
single case shown in Figure~\ref{fig-train-track}, we see that the
singularities of $\tcF$ occur only at the boundary of $W$.  Moreover,
there is only one singular leaf touching $\partial_1W$.  As noted in
Section~\ref{sec-sibling}, there is an involution of $W$ corresponding
to $x \mapsto x, y \mapsto Y$.  This preserves $F$ hence and $\tcF$.  Because of
(\ref{eq-peripheral}), the involution interchanges the cusps via
$(\mu_1\leftrightarrow \lambda_2, \lambda_1 \leftrightarrow -\mu_2)$.  Thus $\tcF$ also has only one
singular leaf touching $\partial_2W$.

Let $\delta_k$ be the curve in $\partial_kW$ corresponding to the boundary of
the singular leaf (in the lamination context, $\delta_k$ is called the
degeneracy slope of $\tcF$).  Because of the involution, we have
\[
\delta_1 = c \mu_1 + d \lambda_1 \quad \mbox{and} \quad \delta_2 = -d \mu_2 + c \lambda_2 \quad \mbox{for some $c,d \in \Z$.}
\]

Returning to the foliation $\cF$ in the fiber surface $\Sigma$, we now
know it has no singularities in the interior.  Moreover, if $\gamma_k$ is
a boundary component of the fiber surface $\Sigma$ lying on $\partial_kW$, the
number of singular prongs is just the intersection number $\abs{\gamma_k
  \cdot \delta_k}$.  From the case of $\phi = 2 x^* + 1$ above, using
Lemma~\ref{lem-filling-slopes} it easily follows that $\delta_1 = \mu_1 $
and $\delta_2 = \lambda_2$.

For a general $\phi = a x^* + b y^*$ over $F$ with $b \geq 0$, by
(\ref{eq-claim}) and Lemma~\ref{lem-filling-slopes}, the fiber surface
has $5^{\epsilon_k(\phi)}$ boundary components along $\partial_k W$ and the
foliation $\cF$ has singularities with $5^{-\epsilon_k(\phi)} (2 a +
(-1)^{k+1} b)$ prongs.  We are now ready to prove
Theorem~\ref{thm-dilatation}.

\begin{proof}[Proof of Theorem~\ref{thm-dilatation}]
  Fix $\phi = a x^* + b y^*$ with $0 \leq b < a$ and $\gcd(a,b) = 1$.  By
  Theorem~\ref{thm-genus}, the filled manifold $W_\phi$ fibers over the
  circle in a unique way, and the fiber genus is $a - 2\epsilon(\phi)$.  An
  easy check shows that if $(a,b) \notin \{ (1,0), (2,1), (3,1) \}$, then
  the genus is at least 2.  Similarly, it is straightforward from the
  above formulae to see that for $(a,b) \notin \{ (2, 1), (3,1), (4,3)
  \}$, there are always at least two prongs at each boundary
  component.  Thus when we Dehn fill to get $W_\phi$, the monodromy $\psi$
  remains pseudo-Anosov with the same invariant foliations and
  dilatation factor.  By \cite{McMullen2000}, that dilatation $\lambda$ is
  the largest real root of
  \[
  \Theta\left(x = t^{\phi(x)}, y = t^{\phi(y)}\right) = t^a + t^{-a} - 1 - t^b - t^{-b}
  \]
  or, multiplying through by $t^a$, of $\Lambda(t) = t^{2a} - t^{a + b} - t^{a} -t^{a-b} + 1$.

  Next, we need to check when the foliations are orientable.  Suppose
  $b$ is odd.  Since only one of $\epsilon_k(\phi) \neq 0$, it follows from the
  above formulae that there are singularities with either $2a+b$ or
  $2a-b$ prongs, and both of these numbers are odd.  Thus the foliation
  is nonorientable.  If $a$ is even, then $b$ is odd and so the
  foliation is again nonorientable.  In the remaining case where $a$ is
  odd and $b$ is even, one has by (\ref{eq-alex-poly}) that the
  Alexander polynomial of $W_\phi$ is related to the polynomial $\Lambda(t)$ defining $\lambda$ by
  \[
  \Delta_{W_\phi} = \Delta_W(x=t^a, y=t^b) = t^{2a} + t^{a + b} + t^{a} +t^{a-b} + 1  = \Lambda(-t)
  \]
  Thus in this case, the spectral radius of the action of $\psi$ on the
  homology of the fiber is the same as the dilatation $\lambda$.  Thus the
  foliation must be orientable (see e.g.~\cite[pg.~5]{LanneauThiffeault2009}).

  Summarizing, for $\phi \notin \{ (1,0), (2, 1), (3,1), (4,3) \}$, the manifold
  $W_\phi$ fibers over the circle with fiber a closed surface of genus at least $2$ and
  pseudo-Anosov monodromy having the claimed data.  By Thurston, the bundle
  $W_\phi$ is thus hyperbolic, completing the proof \cite{ThurstonFibered, Otal96}.
\end{proof}

\section{Minimizing dilatation}

Having shown earlier which surgeries on $W$ minimize volume, we now
turn to determining which minimize the dilatation of the
monodromy, and find:
\begin{theorem} \label{thm-minimizes-dilatation}
  For large $g$, the monodromy of $M_g$ minimizes dilatation among all
  $\Sigma_g$-bundles obtained from filling $W$.
\end{theorem}
\noindent
Moreover, the asymptotic behavior of these minimal dilatations is
\begin{theorem}\label{thm-asymptotic-dilatation}
  Let $\lambda_g$ denote the dilatation of $M_g$.  Then $\lim_{g \to \infty}
  \lambda_{g}^{g} = \frac{3+\sqrt{5}}{2} = 1 + \gamma$, where $\gamma$ is the golden
  ratio.
\end{theorem}

To understand how the dilatation depends on $\phi \in H^1(W; \Z)$, we
follow \cite[\S5]{McMullen2000}.  In particular, for such $\phi$ lying
over the rightmost face of the Thurston norm ball, let $K(\phi)$ be the
dilatation of the monodromy of the corresponding bundle.  If we set
$Z(\phi) = 1/\log(K(\phi)$, then this extends to a nonnegative continuous
function defined on $\setdef{ a x^* + b y^* \in H^1(W; \R)}{ \mbox{$a
    \geq 0$ and $\abs{b} \leq a$}}$.  The function $Z$ is homogeneous with
$Z(c \phi) = c Z(\phi)$ and is strictly concave.   It vanishes exactly along the
rays $b = \pm a$ and is real-analytic away from them.  Thus $Z$ is
completely determined by $z \maps [-1,1] \to \R$ given by $z(t) = Z(x^*
+ t y^*)$.  As noted in Section~\ref{sec-sibling}, there is an
involution of $W$ which acts on $H^1(W;\R)$ by flipping across the
$x^*$-axis.  This involution preserves $K$ and hence $Z$; it follows
that $z$ is even and, since it is strictly concave, has a unique
maximum at $0$.

\begin{proof}[Proof of Theorem~\ref{thm-minimizes-dilatation}]

  Fix a genus $g$, and recall from Section~\ref{sec-examples} that
  the $\phi = a x^* + b y^*$ giving $\Sigma_g$-bundles fall into two
  classes, namely $\A_g$ where $a = g$ and $\E_g$ where $a = g + 2$.
  We seek to minimize $K(\phi)$, or correspondingly maximize $Z(\phi)$.
  From the discussion above, among $\phi$ within just one of $\A_g$ or
  $\E_g$, it is clear that $Z(\phi)$ is maximized when $b$ is least.

  If $g \equiv 3 \bmod 5$, then $\E_g$ is empty and $M_g$ comes from the
  $\phi \in \A_g$ with least $b$.  Thus $Z(\phi)$ is maximized for $M_g$
  and we're done with this case.

  Otherwise, let $\phi_g = (g+2) x^* + b y^*$ denote the element of $\E_g$ with smallest $b$
  so that $M_g = W_{\phi_g}$.  We will show that for large $g$ one has
  \begin{equation}\label{eq-Z-claim}
    Z(g x^*) < Z(\phi_g)
  \end{equation}
  This suffices since we know $Z(\phi) < Z(g x^*)$ for any $\phi \in \A_g$.
  By the homogeneity of $Z$, the claim (\ref{eq-Z-claim}) is
  equivalent to
  \begin{equation}\label{eq-z-claim}
  g z(0) < (g + 2) z\left( \frac{b}{g+2} \right)  \mtext{or}  \frac{g}{g+2} < \frac{z\left(b/(g+2)\right)}{z(0)}
  \end{equation}
  Now $z$ is real-analytic and $z'(0) = 0$ as it has a max at $0$.
  Expanding $z$ in a power series about $0$, and applying the bound on
  $b$ from Proposition~\ref{prop-phi-refined}(\ref{prop-bound-case}),
  we get for large $g$ that
  \[
   \frac{z\left(b/(g+2)\right)}{z(0)}  >  1 - C_1\left( \frac{b}{(g+2)} \right)^2  > 1 - C_2 \frac{\log^2 g}{(g+2)^{2}}.
  \]
  Comparing this with $g/(g + 2) = 1 - 2/g + O(g^{-2})$ gives the
  claim (\ref{eq-z-claim}) for large $g$, proving the theorem.
\end{proof}

\begin{proof}[Proof of Theorem~\ref{thm-asymptotic-dilatation}]
  Let $\phi_g$ be the class giving $M_g$.  The properties of $Z$
  mentioned above mean that $K$ also extends to a continuous function
  with $K( c \phi ) = K(\phi)^{1/c}$.  Thus
  \[
  \lambda_g^g = K(\phi_g)^g = K\left({\scriptstyle \frac{1}{g}} \phi_g\right)
  \]
  By Proposition~\ref{prop-phi-refined}(\ref{prop-bound-case}), we
  know that if $\phi_g = a x^* + b y^*$ then $b$ is $O(\log g)$.  Thus
  we see that $(1/g)\phi_g$ converges to $x^*$ as $g \to \infty$.  Hence by continuity
  $\lim_{g \to \infty} \lambda_g^g = K(x^*)$.  The latter is the largest root of
  \[
  \Theta(x = t, y = 1) = t^2 - 3 t + 1
  \]
  which is $(3 + \sqrt 5)/2 = 1 + \gamma$ as claimed.
\end{proof}

We now turn to our partial answer to Question~\ref{quest-simult-min}.

\begin{theorem}\label{thm-min-vol-dilatation-differ}
  Assume Conjecture~\ref{conj-closed-b12}.  Let $\delta_g$ be the smallest
  dilatation of any pseudo-Anosov of $\Sigma_g$.  Then for all large $g \equiv
  3 \bmod 5$, no pseudo-Anosov $\psi_g$ with dilatation $\delta_g$ gives a
  hyperbolic $\Sigma_g$-bundle of least volume.
\end{theorem}

\begin{proof}
  Assuming Conjecture~\ref{conj-closed-b12} gives by
  Theorem~\ref{thm-main-vol} that we know there is a unique
  $\Sigma_g$-bundle of least volume when $g$ is large, namely $M_g$.  Let $\lambda_g$ be the
  corresponding dilatation.  We compare this to the pseudo-Anosov of
  $\Sigma_g$ constructed by Hironaka in \cite[Thm.~1.4]{Hironaka2009},
  whose dilatation we denote by $\epsilon_g$.  In particular, we just need
  to show that $\epsilon_g < \lambda_g$.

  Following the notation of \cite{Hironaka2009}, let $\lambda_{(b,a)}$ be
  the largest root of the polynomial in Theorem~\ref{thm-dilatation}.
  Then the dilatation $\epsilon_g$ is either $\lambda_{(1, g+1)}$ or $\lambda_{(3,
    g+1)}$.  As $g \equiv 3 \bmod 5$, Corollary~\ref{cor-dilatations}
  gives that $\lambda_g = \lambda_{(b, g)}$ where $b \in \{2, 3, 5\}$.  By
  Propositions~4.1 and 4.2 of \cite{Hironaka2009}, one has
  \[
  \epsilon_g \leq \lambda_{(3, g+1)} \leq \lambda_{(1, g)} < \lambda_g.
  \]
  as needed.
\end{proof}

Finally, we identify the minimal dilatation $\delta_g^+$ of any
pseudo-Anosov of $\Sigma_7$ with orientable invariant foliations.

\begin{proof}[Proof of Theorem~\ref{thm-delta-plus-7}]
  Consider the monodromy $\psi$ of the bundle $M_7$, which corresponds to $\phi
  = 9 x^* + 2 y^*$.  By Theorem~\ref{thm-dilatation}, the invariant
  foliations of $\psi$ are orientable, and the dilatation $\lambda_7 \approx
  1.11548110945659$ is the largest root of
  \begin{multline*}
  t^{18} - t^{11} - t^9 - t^7 + 1 = \left(t^{4} - t^{3} + t^{2} - t + 1\right) \cdot \\ 
  \left(t^{14} + t^{13} - t^{9} - t^{8} - t^{7} - t^{6} - t^{5} + t + 1\right)
  \end{multline*}
  Now $\lambda_7$ is a root of the second irreducible factor above,
  which is exactly the polynomial used in \cite{LanneauThiffeault2009}
  to give a lower bound on $\delta_7^+$.   Thus $\delta_7^+ = \lambda_7$, as claimed.
\end{proof}

{\RaggedRight \bibliographystyle{math-1.4} \bibliography{small-bundles.bib} }

\end{document}